\documentclass[a4paper,10pt]{article}

\def\u{\bigsqcup}
\def\eps{\varepsilon}
\title{ Simple Spectrum for Tensor  Products of  Mixing Map Powers  }
\author{V.V. Ryzhikov\footnote{\large This work is partially  supported by  the grant
NSh 8508.2010.1.}}

\begin{document}
\Large
\maketitle

\section{Introduction} 
In this note we  consider  measure-preserving transformations  of a Probability space $(X,\mu$).
We  prove the existence of a mixing rank one construction $T$ such that     the product
 $T\otimes T^2\otimes T^3\otimes\dots$ has simple spectrum. This result has been  announced in \cite{R}. It had  an application in recent Thikhonov's proof \cite{T} of the existence of mixing transformation with homogeneous spectrum of  multiplicity   $m>2$ (see \cite{D}).  
Let us remark that for generic non-mixing  transformations the above spectral properties    have been found  by Ageev  \cite{Ag}.

{\bf Rank one construction} is determined by $h_1$ and a sequence $r_j$  of cuttings and  a sequence $\bar s_j$ of spacers
$$ \bar s_j=(s_j(1), s_j(2),\dots, s_j(r_j-1),s_j(r_j)).$$  We recall its definition.
Let our  $T$ on the step $j$  be associated with  a collection of disjoint sets (intervals)
$$E_j, TE_j T^2,E_j,\dots, T^{h_j}E_j.$$
We cut $E_j$ into $r_j$ sets (subintervals)  of the same measure
$$E_j=E_j^1\u E_j^2\u  E_j^3\u\dots\u E_j^{r_j},$$  
then for all $i=1,2,\dots, r_j$ we  consider columns
$$E_j^i, TE_j^i ,T^2 E_j^i,\dots, T^{h_j}E_j^i.$$
Adding $s_j(i)$ spacers we obtain  
$$E_j^i, TE_j^i T^2 E_j^i,\dots, T^{h_j}E_j^i, T^{h_j+1}E_j^i, T^{h_j+2}E_j^i, \dots, T^{h_j+s_j(i)}E_j^i$$
(the above intervals   are disjoint).
For all  $i<r_j$ we set
$$TT^{h_j+s_j(i)}E_j^i = E_j^{i+1}.$$ 
Now we obtain a tower 
$$E_{j+1}, TE_{j+1} T^2 E_{j+1},\dots, T^{h_{j+1}}E_{j+1},$$
where 
 $$E_{j+1}= E^1_j,$$
$$T^{h_{j+1}}E_{j+1}=T^{h_j+s_j(r_j)}E_j^{r_j}, $$
$$h_{j+1}+1=(h_j+1)r_j +\sum_{i=1}^{r_j}s_j(i).$$
So step by step we define a general rank one construction.

\bf On notation. \rm We denote weak operator approximations by $\approx_w$ and $\approx_s$ for strong ones.
$\Theta$ is the orthogonal projection into the space of constant functions in $L_2(X,\mu)$.
Thus, the  expression $T^m\approx_w\Theta$ (for  large $m$) means that $T$ is mixing.
 
\bf Stochastic Ornstein's rank one transformation. \rm D. Ornstein has proved the mixing for almost all spesial stochastic constructions. His proof can be very  shortly presented  in the  following manner.
Let $H_j\to \infty$,  $H_j<<r_j$,  we consider  uniformly distributed stochastic variables
$a_j(i)\in \{0, 1, \dots, H_j\}$ and let 
$$s_j(i)=H_j+a_j(i)- a_j(i+1).$$

Then for $m\in [h_j,h_{j+1})$
$$T^m\approx_w  D_1T^{k_1}P_1 +D_2T^{k_2}P_2 +D_3T^{k_3}P_3,$$ 
where $D_i$ are operators of multiplication by indicators of special parts of $j$-towers
(all $D_i$ and $P_i$ depend on $m$), $|k_1|<h_{j+1},$ $|k_2|,|k_3|<h_j$,  the operators $P_i$ have the form:
$$P_1=\sum_{n\in [-H_{j+1}, H_{j+1}]}c_{j+1}(n)T^{ n}, \ \  c_{j+1}(n)=\frac{H_{j+1}+1-|n|}{(H_{j+1}+1)^2}, $$
$$P_{2,3}=\sum_{n\in [-H_{j}, H_{j}]}c_{j}(n)T^{ n}, \ \  c_{j}(n)=\frac{H_{j}+1-|n|}{(H_{j}+1)^2}. $$
 We have $$\|D_iT^{k_i}\|\leq 1$$ and $$P_i\approx_s \Theta$$
since $T$ is ergodic.  Finally we get for large $m$ 
$$T^m\approx_w\Theta.$$
\section{ $T^{sm_j}\to\Theta $ for  a given     exclusive $s$}
  Weak limits  $aI +(1-a)\Theta$ are well known in ergodic theory. They have been   used   in connection with 
 Kolmogorov's problem  \cite{S}and   for a machinery of counterexamples \cite{JL}. 

\vspace{2mm}
\bf LEMMA 1. \it  For any $\eps>0$,  $N$ and  $s\in[1, N]$ there is a rank one $(1-\eps)$-partially mixing  construction $T$ with the following property:  for a sequence $m_j$ 
$$ {\bf (N,s)-Property}  \left\{ 
\begin{array}{rcl}
T^{sm_j}\to   \Theta, \\   
T^{km_j}\to  (1-a_k)\Theta  + a_k I\\
\end{array}
\right. 
 $$  
for some  $a_k>0$, $k\neq s,$   $1\leq k\leq N$.
\rm
\vspace{3mm}
\vspace{3mm}

We are able to work with staircase spacer arrays \cite{A} as well as algebraic spacers \cite{R1}, but we prefer  stochastic constructions \cite{O}. We do not try to construct  explicit examples here and  follow this simple way.
 
Proof. Let $s=3$, $N=5$. A sequence of spacers is organized as follows.  A spacer vector $\bar{s_j}$ is
a concatenation of arrays 
$$S1,S1, A1,S2,S2,A2, S4,S4,A4,S5,S5,A5,$$
where $Sk,Ak$ are independent arrays of spacers. 
Moreover let
 arrays $Sk$  be  stochastic Ornstein's spacer sequences of the length $kL_j$ with an  average value equals to $H_j$; let an array $Ak$  be of a length $[\eps^{-1}kL_j]$.

Let $m_j=(h_j+H_j)L_j$, then for a small constant $a>0$ (we omit its calculation) and 
$k\neq s$,   $1\leq k\leq N$,  one gets  
$$T^{km_j}\approx_w  ka_kI +(1-ka_k)( D_1T^{k_1}P_1 +D_2T^{k_2}P_2 +D_3T^{k_3}P_3)\approx_w  $$
$$\approx_w  ka_kI +(1-ka_k)\Theta,  \ \ a_k>a.  $$
Via Ornstein's approach   the mixing  is everywhere   in $j$-tower except  a part $D$  that is situated  under the second spacer array $Sk$.  For this part we have  for measurable sets $B,B'$ 
$$\mu(T^{km_j}B\cap B'\cap D)\approx \mu(D)\mu(B\cap B'),$$
so   $ ka_kI $ appears.
However 
$$T^{3m_j}\approx_w\Theta$$
since we "forget" to  copy an  array of the length $3L$. 
\section{ Exclusive $n$ for which  $T^{nm_j}\to aI+bT+c\Theta$}

{\bf $(n,a,b)$-constructions.}  Let $r_j\to\infty$.  We fix positive $a,b$,  $a+b+c=1$, and $n>1$. 
For a subsequence $r_{j'-1}$  we produce  a flat part (a-part), a polynomial part (b-part) and a  mixing  part (c-part) (stochastic \cite{O}, algebraic \cite{R1}, or staircase \cite{A} that we use here). These parts will be now  provided by the following  spacer sequence 
 $\bar s_{j'}$.  We set ( writing  again $j$ instead of $j'$)

{\bf a-Part:} for  $i=1, 2, \dots, [a r_j]$ 
$$s_j(i) = H_j.$$

{\bf b-Part:}  for  $i\in  ([a r_j],  [(a+b) r_j])$ if  $i=ni'$ we set
$s_j(i)=nH_j-1,$
otherwise
$s_j(i)=0.$
So this part of  spacer vector  looks as 

$\dots,0,0,\dots,0, nH_j-1,\ 0,0,\dots,0, nH_j-1,\ 0,0,\dots,0, nH_j-1,\ 0, \dots$

{\bf Mixing c-Part:} $s_j(i)=i$ for $i>[(a+b) r_j]$.

{\bf A condition for $(j-1)$-steps.} \it We define  on $j-1$-step our construction to be  a pure staircase and we set $H_j=h_{j-1}$
(recall that $j=j'$ is a subsequence).   \rm
\vspace{3mm}

 {\bf Weak limits. } Let $m_j=h_j+H_j$. We get the following convergences:
$${\bf  n-Property}\ \ \ \ \ \left\{
\begin{array}{rcl}
T^{m_j}\to   aI+(b  +  c)\Theta,\\
T^{2m_j}\to  aI+(b  +  c)\Theta,\\
\dots \\
T^{(n-1)m_j}\to  aI+(b  +  c)\Theta,\\
T^{nm_j}\to  aI+ b T + c \Theta.\\
\end{array}
\right.
$$
Indeed, we have
$$T^{Km_j}\approx_w  aT^0+  b\left(\frac{n-K}{n}T^{KH_j} + \frac{K}{n}T^{(K-n)H_j+1}\right)  
 +\frac{1}{r_j}\sum_{i>(a+b)r_j}^{r_j} T^{-2i-1},  
$$
$$T^{Km_j}\ \approx_w \ \ 
aI  +b\left(\frac{n-K}{n}T^{KH_j} + \frac{K}{n}T^{(K-n)H_j+1}\right)  + c\Theta.$$
For $K=n$
$$b\left(\frac{n-K}{n}T^{KH_j} + \frac{K}{n}T^{(K-n)H_j+1}\right) =bT.$$
For  $K=1,2,\dots, n-1$ we use  $T^{KH_j},T^{KH_j+1}\approx_w\Theta$ and obtain  
$$b\left(\frac{n-K}{n}T^{KH_j} + \frac{K}{n}T^{(K-n)H_j+1}\right)\approx_w \ b\Theta.$$

\section{Main result}

\bf LEMMA 2. \it   Let  for  $m=2,3,\dots, n$ and all $s\leq m$ a transformation $T$ have  $m$-Properies and  $(s,m)$-Properties. Then 
$T\otimes T^2\otimes\dots T^n$ has simple spectrum.\rm
\vspace{3mm}

Proof. A cyclic vector for $T$ in $H=L^0_2$ is denoted by $f$.
We shall prove that a cyclic space $C_F$ is $H^{\otimes n}$, where  $F=f^{\otimes n}$ and  $T\otimes T^2\otimes\dots T^n$ is restricted to  $H^{\otimes n}$.  
For $S=T\otimes T^2\otimes\dots T^{n-1}$ we assume
it has  simple spectrum by induction.  

From $n$-property  we get    
$$b^{n-1}f^{\otimes n-1}\otimes (aI+bT)f \in C_F,$$
hence,
$f^{\otimes n-1}\otimes Tf \in C_F,$
thus, for all $k$
$$f^{\otimes n-1}\otimes T^kf \in C_F.$$
This implies
$$f^{\otimes n-1}\otimes H \subset C_F, \ \
S^if^{\otimes n-1}\otimes H \subset C_F,
\ \ H^{\otimes n-1}\otimes H \subset C_F.$$

To see that $T\otimes T^2\otimes\dots T^n$ has a simple spectrum in $L_2^{\otimes n}$ we note
that all different products $T^{n_1}\otimes \dots \otimes T^{n_k}$ are spectrally disjoint.
This follows directly from (s,n)-Properties (see Lemma 1).
For example,  if $s=2$ (in Lemma 2), then   in $H^{\otimes 3}$
$$(T\otimes T^2\otimes\ T^5)^{m_j}\ \to_w\  0,$$ but
$$(T\otimes T^3\otimes\ T^5)^{m_j}\ \to_w \ a_1a_3a_5I>a^3I.$$
\vspace{3mm}

{\bf THEOREM.} \it  There is a mixing transformation $T$ such that $T\otimes T^2\otimes\ T^3 \dots$
has  simple spectrum. \rm
\vspace{3mm}

Proof.  We construct  rank one transformations  $T_p$ with  n-Properties and $(s,N)$-Properties ($n,N\leq p$). We make  these constructions $c_p$-partially
mixing with  $T_p\otimes T_p^2\otimes\ T_p^3 \dots$  of  simple spectrum (Lemma 2).
Then  $c_p$  tends very slowly to 1, and we  force a  limit mixing construction $T$ to have the  desired  spectral property  via standard technique (see \cite{D}, \cite{R}).   In \cite{R} we define
 $T_p\to T$ ($p\to\infty$) to have simple spectrum for all $T^{\odot n}$.  Replacing this aim
by another one we provide  simple spectrum of $T\otimes T^2\otimes \dots T^n$  via  the same methods.

Finally let us formulate   a similar  problem on flows:

\bf Conjecture.  \it  There is a mixing flow $T_t$ such that for all collections of  different $t_i>0$  the  products
$T_{t_1}\otimes T_{t_2}\otimes T_{t_3}\dots$
have  simple spectrum. Moreover the same is true for 
$$exp(T_{t_1})\otimes exp(T_{t_2})\otimes exp(T_{t_3})\dots$$
(here  $T_{t_i}$ are now treated  as unitary operators  restricted  onto $H$). 
\rm

The main difficulty is not to find a solution  but is to find  an  elegant one.     It seems that the following lemma could be useful.
\vspace{2mm}

\bf LEMMA 3. \it If for a flow $T_t$ with simple spectrum   and  any positive different  $s, t_1,t_2,\dots, t_n$ there is $m_j\to \infty$ and positive $a_1,a_2,\dots, a_n$  such that  
$$ 
T_{t_im_j}\to a_iI +(1-a_i)  \Theta, \  i=1,2,\dots, n-1,$$  
$$T_{t_nm_j}\to  a_nT_s + (1-a_n)\Theta,  $$
then the corresponding (Gaussian) automorphisms 
$$exp(T_{t_1})\otimes exp(T_{t_2})\otimes exp(T_{t_3})\dots$$ 
have  simple spectrum. 
\rm

Remark.  There is a weakly mixing flow $T_t$ possessing the following property: given  $a\in [0,1]$ there is a  sequence 
$m_i$ such that  for any real $s>0$ there is a subsequence $m_{i(k)}$ (it depends on $s$) providing 
$$T_{sm_{i(k)}}\to aI+(1-a)\Theta.$$
(It is not possible to have $T_{sm_{i}}\to aI+(1-a)\Theta$ for a  set of $s$ of a positive measure.)\large

Hint: let us consider  rank one flows with $r_j>>>h_j$.

 E-mail: vryz@mail.ru
\end{document}